\documentclass[a4paper,11pt]{amsart}
\usepackage{amsfonts,amsmath,amsthm,amssymb,multicol}
\usepackage{enumitem}
\usepackage{mathrsfs}
\usepackage[margin=1.15in]{geometry}
\usepackage[colorlinks=true,linkcolor=blue,filecolor=blue,citecolor=red]{hyperref}
\usepackage{orcidlink}
\usepackage{cleveref}

\theoremstyle{plain}
\newtheorem{theorem}{Theorem}[section]

\theoremstyle{definition}

\begin{document}
\title[Congruences for the difference between even cranks and odd cranks]{A note on congruences for the difference between even cranks and odd cranks}
\author[Russelle Guadalupe]{Russelle Guadalupe\orcidlink{0009-0001-8974-4502}}
\address{Institute of Mathematics, University of the Philippines Diliman, Quezon City 1101, Philippines}
\email{rguadalupe@math.upd.edu.ph}

\renewcommand{\thefootnote}{}

\footnote{2020 \emph{Mathematics Subject Classification}: Primary 11P83, Secondary 05A17, 11P81}

\footnote{\emph{Key words and phrases}: crank of a partition, Ramanujan-type congruences, $q$-series, dissection formulas}

\renewcommand{\thefootnote}{\arabic{footnote}}

\setcounter{footnote}{0}

\begin{abstract}
Recently, Amdeberhan and Merca proved some arithmetic properties of the crank parity function $C(n)$ defined as the difference between the number of partitions of $n$ with even cranks and those with odd cranks and the sequence $a(n)$ whose generating function is the reciprocal of that of $C(n)$. The function $C(n)$ was first studied by Choi, Kang, and Lovejoy. In this note, we give new elementary proofs of some of their main results and extend them. In particular, we establish Ramanujan-type congruences modulo $5$ and $25$ for certain finite sums involving $C(n)$ and $a(n)$.  Our proofs employ the results of Cooper, Hirschhorn, and Lewis, and certain identities involving the Rogers-Ramanujan continued fraction $R(q)$ due to Chern and Tang.
\end{abstract}

\maketitle

\section{Introduction} \label{sec1}

Throughout this paper, let $f_m := \prod_{n\geq 1}(1-q^{mn})$ for a positive integer $m$ and a complex number $q$ with $|q| < 1$. Recall that a partition of a positive integer $n$ is a finite non-increasing sequence of positive integers, known as its parts, whose sum is $n$. It is well-known that the generating for the number $p(n)$ of partitions of $n$ with $p(0):=1$ is given by 
\begin{align*}
\sum_{n\geq 0} p(n)q^n = \dfrac{1}{f_1}.
\end{align*}
In his 1919 paper, Ramanujan \cite{rama} discovered three congruences satisfied by $p(n)$, namely,
\begin{align*}
p(5n+4)&\equiv 0\pmod{5},\\
p(7n+5)&\equiv 0\pmod{7},\\
p(11n+6)&\equiv 0\pmod{11}.
\end{align*}
Ramanujan \cite{rama} and Darling \cite{dar} independently proved the first two of these congruences, which are now special cases of the following remarkable result due to Watson \cite{wat} and Atkin \cite{atk} that reads
\begin{align*}
p(n)\equiv 0\pmod{\ell^\beta}
\end{align*}
for all $\alpha\geq 1$ and $n\geq 1$ such that $24n\equiv 1\pmod{\ell^\alpha}$, where $\beta := \alpha$ if $\ell\in \{5,11\}$ and $\beta := \lfloor\alpha/2\rfloor+1$ if $\ell=7$.

Stemming from the previous works of Dyson \cite{dys} and Garvan \cite{garv}, Andrews and Garvan \cite{andg} introduced in 1988 the crank of a partition to provide a unified combinatorial interpretation of the above congruences of Ramanujan. Given a partition $\lambda$, let $\omega(\lambda)$ be the number of $1$'s in $\lambda$, $l(\lambda)$ be the largest part of $\lambda$, and $\mu(\lambda)$ be the number of parts of $\lambda$ larger than $\omega(\lambda)$. Then the crank of $\lambda$ is defined by 
\[c(\lambda) = \begin{cases}
l(\lambda)  & \text{if }\omega(\lambda) = 0,\\
\mu(\lambda)-\omega(\lambda)  & \text{if }\omega(\lambda) > 0.
\end{cases}\]

In 2009, Choi, Kang, and Lovejoy \cite{choi} introduced the crank parity function $C(n)$ defined as the difference between the number of partitions of $n$ with even cranks and those with odd cranks and showed that the generating function of $C(n)$ is given by 
\begin{align*}
\sum_{n\geq 0} C(n)q^n = \dfrac{f_1^3}{f_2^2}.
\end{align*}
Andrews \cite{and} independently discovered the above generating function of $C(n)$. Using the theory of modular forms, Choi, Kang, and Lovejoy \cite{choi} established an infinite family of congruences modulo powers of $5$ for $C(n)$ and the generating function of $C(5n+4)$. 

\begin{theorem}{\cite[Theorem 1.1]{choi}}\label{thm11}
For all $\alpha \geq 0$, we have $C(n)\equiv 0\pmod{5^{\alpha+1}}$ if $24n\equiv 1\pmod{5^{2\alpha+1}}$. 
\end{theorem}

\begin{theorem}{\cite[Theorem 1.2]{choi}}\label{thm12}
We have the identity
\begin{align*}
\sum_{n\geq 0} C(5n+4)q^n = 5\dfrac{f_1^2f_5f_{10}^2}{f_2^4}.
\end{align*}
\end{theorem}

Recently, Tang \cite{tang} gave elementary proofs of Theorems \ref{thm11} and \ref{thm12} and found more congruences for $C(n)$ modulo small powers of $5$. Very recently, Amdeberhan and Merca \cite{amdb} studied $C(n)$ and another sequence $a(n)$ whose generating function is the reciprocal of that of $C(n)$:
\begin{align*}
\sum_{n\geq 0} a(n)q^n = \dfrac{f_2^2}{f_1^3}.
\end{align*}
Combinatorically, $a(n)$ is the number of partitions of $n$ where each odd part is marked with three distinct colors. Utilizing certain $5$-dissection formulas for Ramanujan's theta functions, Amdeberhan and Merca \cite{amdb} obtained another elementary proof of Theorem \ref{thm12}. They also obtained another congruence modulo $5$ for a finite sum involving $C(n)$ and a congruence modulo $7$ for $a(n)$ using the \textit{Mathematica} package \texttt{RaduRK} devised by Smoot \cite{smoot}, which implements the Ramanujan-Kolberg algorithm found by Radu \cite{radu} via the theory of modular forms. 

\begin{theorem}{\cite[Corollary 18(b)]{amdb}}\label{thm13}
For all $n\geq 0$, we have 
\begin{align*}
\dfrac{1}{5}\sum_{k=-\infty}^\infty (1+6k)C\left(50n+49-\dfrac{25k(3k+1)}{2}\right)\equiv 0\pmod{5}.
\end{align*}
\end{theorem}

\begin{theorem}{\cite[Theorem 6]{amdb}}\label{thm14}
For all $n\geq 0$, we have $a(7n+2)\equiv 0\pmod{7}$.
\end{theorem}

The goal of this note is to provide elementary proofs of Theorems \ref{thm13} and  \ref{thm14} and extend them. In particular, we derive congruences modulo $5$ and $25$ for certain finite sums involving $C(n)$ and $a(n)$. 

\begin{theorem}\label{thm15}
For all $n\geq 0$, we have
\begin{align}
\sum_{k\geq 0} a\left(25n+16-\dfrac{5k(k+1)}{2}\right)&\equiv 0\pmod{5},\label{eq11}\\
\dfrac{1}{5}\sum_{k\geq 0} C\left(125n+114-\dfrac{5k(k+1)}{2}\right)&\equiv 0\pmod{25}.\label{eq12}
\end{align}
\end{theorem}

\begin{theorem}\label{thm16}
Let $p\equiv 13, 17,19, 23\pmod{24}$ be a prime. Then for all $n\geq 0$ and $r\in \{1,\ldots, p-1\}$, we have
\begin{align*}
\sum_{k=-\infty}^\infty (-1)^ka\left(5p^2n+5pr+\dfrac{25p^2-1}{24}-5k^2\right)&\equiv 0\pmod{5}.
\end{align*}
\end{theorem}

We organize the rest of the paper as follows. We first prove in Section \ref{sec2} Theorem \ref{thm14} using a special case of the result of Cooper, Hirschhorn, and Lewis \cite{cophir}. We next recall in Section \ref{sec3} important formulas, including the $5$-dissections of $f_1$ and $1/f_1$ in terms of the Rogers-Ramanujan continued fraction defined by
\begin{align*}
	R(q) = \prod_{n\geq 1}\dfrac{(1-q^{5n-4})(1-q^{5n-1})}{(1-q^{5n-3})(1-q^{5n-2})}
\end{align*} 
(here, we omit the extra factor $q^{1/5}$ in the definition of $R(q)$ for convenience) and certain identities involving $R(q)$ and the parameter
\begin{align*}
K := \dfrac{f_2f_5^5}{qf_1f_{10}^5}
\end{align*}
due to Chern and Hirschhorn \cite{chehir} and Chern and Tang \cite{chetn}. We then apply these identities to demonstrate Theorems \ref{thm13} and \ref{thm15} in Section \ref{sec4}, and, together with \cite{cophir}, Theorem \ref{thm16} in Section \ref{sec5}. We finally give some remarks in Section \ref{sec6}.

\section{Proof of Theorem \ref{thm14}}\label{sec2}

\begin{proof}[Proof of Theorem \ref{thm14}]
Let $f_1^4f_2^2:=\sum_{n\geq 0} d(n)q^n$ so that
\begin{align}
	\sum_{n\geq 0} a(n)q^n = \dfrac{f_2^2}{f_1^3}\equiv \dfrac{f_1^4f_2^2}{f_7}\equiv \dfrac{1}{f_7}\sum_{n\geq 0} d(n)q^n \pmod{7}\label{eq21}
\end{align}
by the binomial theorem. We consider the terms of (\ref{eq21}) involving $q^{7n+2}$, divide both sides by $q^2$, and then replace $q^7$ with $q$. We obtain
\begin{align}
	\sum_{n\geq 0} a(7n+2)q^n \equiv \dfrac{1}{f_1}\sum_{n\geq 0} d(7n+2)q^n \pmod{7}.\label{eq22}
\end{align}
Observe that $d(2) = d(9) = 0$ and 
\begin{align}
 d(7n+16) = 49d\left(\dfrac{n}{7}\right)\label{eq23}
\end{align}
for all $n\geq 0$ by \cite[Theorem 2]{cophir}, where we denote $d(n/7) = 0$ if $7\nmid n$. Combining (\ref{eq22}) and (\ref{eq23}) leads to 
\begin{align*}
	\sum_{n\geq 0} a(7n+2)q^n \equiv \dfrac{1}{f_1}\sum_{n\geq 0} d(7n+16)q^{n+2}\equiv 0 \pmod{7},
\end{align*}
so that $a(7n+2)\equiv 0\pmod{7}$ as desired.
\end{proof}

\section{Identities involving $R(q)$ and $K$}\label{sec3}

In this section, we require important identities needed prove the remaining theorems. We begin with the following $5$-dissections \cite[(8.1.4), (8.4.4)]{hirsc}
\begin{align}
	f_1 &= f_{25}\left(\dfrac{1}{R_5}-q-q^2R_5\right),\label{eq31}\\
	\dfrac{1}{f_1} &= \dfrac{f_{25}^5}{f_5^6}\left(\dfrac{1}{R_5^4}+\dfrac{q}{R_5^3}+\dfrac{2q^2}{R_5^2}+\dfrac{3q^3}{R_5}+5q^4-3q^5R_5+2q^6R_5^2-q^7R_5^3+q^8R_5^4\right)\label{eq32},
\end{align}
where $R_m := R(q^m)$, and the following identities \cite[(9.10), (9.11)]{chehir}
\begin{align}
	K+1 &= \dfrac{f_2^4f_5^2}{qf_1^2f_{10}^4},\label{eq33}\\
	K-4 &= \dfrac{f_1^3f_5}{qf_2f_{10}^3}.\label{eq34}
\end{align}
Let 
\begin{align*}
P(m,n) := \dfrac{1}{q^mR_1^{m+2n}R_2^{2m-n}} + (-1)^{m+n}q^mR_1^{m+2n}R_2^{2m-n}
\end{align*}
for any $m\in\mathbb{N}$ and $n\in\mathbb{Z}$. 
Chern and Tang \cite[Theorem 1.1]{chetn} established the recurrence relations
\begin{align}
	P(m,n+1) &= 4K^{-1}P(m,n) + P(m,n-1),\label{eq35}\\
	P(m+2,n) &= KP(m+1,n)+P(m,n),\label{eq36}
\end{align}
with the initial values
\begin{align}
	P(0,0) &= 2,\label{eq37}\\
	P(0,1) &= 4K^{-1},\label{eq38}\\
	P(1,0) &= K,\label{eq39}\\
	P(1,-1) &= 4K^{-1}-2+K.\label{eq310}
\end{align}
We finally note that $f_m^{5^k} \equiv f_{5m}^{5^{k-1}}\pmod{5^k}$ for all $k\geq 1$ and $m\geq 1$, which follows from the binomial theorem and will be frequently used without comment.

\section{Proofs of Theorems \ref{thm13} and \ref{thm15}}\label{sec4}

\begin{proof}[Proof of Theorem \ref{thm13}]
Let
\begin{align*}
	\sum_{n\geq 0} A(n)q^n := \dfrac{f_1^2f_5^6}{f_2^4}.
\end{align*}
In view of Theorem \ref{thm12} and the identity of Ramanujan \cite[(10.7.3)]{hirsc}
\begin{align*}
	\sum_{k =-\infty}^\infty (6k+1)q^{k(3k+1)/2} = \dfrac{f_1^5}{f_2^2},
\end{align*}
we have
\begin{align*}
	A(n) = \dfrac{1}{5}\sum_{k =-\infty}^\infty (1+6k)C\left(5n+4-\dfrac{25k(3k+1)}{2}\right)
\end{align*}
so it suffices to prove that $A(10n+9)\equiv 0\pmod{5}$ for all $n\geq 0$. Applying (\ref{eq31}) and (\ref{eq32}) on the generating function of $A(n)$ yields, modulo $5$,
\begin{align}
\sum_{n\geq 0} A(n)q^n \equiv \dfrac{f_2f_1^2f_5^6}{f_{10}}\equiv\dfrac{f_5^6}{f_{10}}\cdot f_{25}^2f_{50}\left(\dfrac{1}{R_5}-q-q^2R_5\right)^2\left(\dfrac{1}{R_{10}}-q^2-q^4R_{10}\right).\label{eq51}
\end{align}
We extract the terms of (\ref{eq51}) involving $q^{5n+4}$, divide both sides by $q^4$, and then replace $q^5$ with $q$. We infer from (\ref{eq34}) and (\ref{eq35}) that
\begin{align*}
\sum_{n\geq 0} A(5n+4)q^n&\equiv \dfrac{f_1^6f_5^2f_{10}}{f_2}(1+P(0,-1))\equiv \dfrac{f_1f_5^3f_{10}}{f_2}\cdot\dfrac{K-4}{K}\\
&\equiv  \dfrac{f_1f_5^3f_{10}}{f_2}\dfrac{f_1^3f_5}{qf_2f_{10}^3}\cdot \dfrac{qf_1f_{10}^5}{f_2f_5^5}\equiv f_2^2f_{10}^2\pmod{5}.
\end{align*}
Since the $q$-expansion of $f_2^2f_{10}^2$ contains only terms with even exponents, we deduce that
\begin{align*}
	\sum_{n\geq 0} A(10n+9)q^n \equiv \sum_{n\geq 0} A(5(2n+1)+4)q^n \equiv 0\pmod{5},
\end{align*}
which yields the desired congruence.
\end{proof}

\begin{proof}[Proof of Theorem \ref{thm15}]
We first show (\ref{eq11}) by finding the generating function modulo $5$ of $a(5n+1)$. From (\ref{eq31}), we have that
\begin{align}
\sum_{n\geq 0} a(n)q^n \equiv \dfrac{f_1^2f_2^2}{f_5}\equiv \dfrac{f_{25}^2f_{50}^2}{f_5}\left(\dfrac{1}{R_5}-q-q^2R_5\right)^2\left(\dfrac{1}{R_{10}}-q^2-q^4R_{10}\right)^2\pmod{5}\label{eq52}.
\end{align}
We extract the terms of (\ref{eq52}) involving $q^{5n+1}$, divide both sides by $q$, and then replace $q^5$ with $q$. We see from (\ref{eq33}), (\ref{eq35}), (\ref{eq38}), and (\ref{eq39}) that 
\begin{align*}
\sum_{n\geq 0} a(5n+1)q^n &\equiv q\dfrac{f_5^2f_{10}^2}{f_1}(1-2P(0,-1)-2P(1,0))\equiv q\dfrac{f_5^2f_{10}^2}{f_1}\left(1+\dfrac{8}{K}-2K\right)\\
&\equiv 3q\dfrac{f_5^2f_{10}^2}{f_1}\cdot\dfrac{(K+1)^2}{K}\\
&\equiv 3q\dfrac{f_5^2f_{10}^2}{f_1}\left(\dfrac{f_2^4f_5^2}{qf_1^2f_{10}^4}\right)^2\left(\dfrac{qf_1f_{10}^5}{f_2f_5^5}\right)\equiv 3\dfrac{f_2^7f_5}{f_1^4f_{10}}\equiv 3f_1f_2^2\pmod{5}.
\end{align*}
We recall the identity \cite[(1.5.3)]{hirsc}
\begin{align*}
\dfrac{f_2^2}{f_1} = \sum_{k\geq 0} q^{k(k+1)/2}
\end{align*}
so that
\begin{align}
\dfrac{3f_{10}}{f_2} \equiv \dfrac{f_2^2}{f_1}\cdot 3f_1f_2^2&\equiv \left(\sum_{k\geq 0} q^{k(k+1)/2}\right)\left(\sum_{j\geq 0} a(5j+1)q^j\right)\nonumber\\
&\equiv \sum_{n\geq 0}\left(\sum_{k\geq 0}a\left(5n+1-\dfrac{5k(k+1)}{2}\right)\right)q^n\pmod{5}.\label{eq53}
\end{align}
We apply (\ref{eq32}) on (\ref{eq53}) and look at the terms of the resulting expression involving $q^{5n+3}$. We arrive at
\begin{align*}
\sum_{n\geq 0}\left(\sum_{k\geq 0}a\left(25n+16-\dfrac{5k(k+1)}{2}\right)\right)q^n\equiv 15q\dfrac{f_{10}^5}{f_2^5}\equiv 0\pmod{5},
\end{align*}
yielding (\ref{eq11}). We next prove (\ref{eq12}); note that by Theorem \ref{thm12},
\begin{align}
\dfrac{f_1f_5f_{10}^2}{f_2^2} = \dfrac{f_2^2}{f_1}\cdot\dfrac{f_1^2f_5f_{10}^2}{f_2^4}&=\left(\sum_{k\geq 0} q^{k(k+1)/2}\right)\left(\dfrac{1}{5}\sum_{j\geq 0} C(5j+4)q^j\right)\nonumber\\
&=\sum_{n\geq 0}\left(\dfrac{1}{5}\sum_{k\geq 0}C\left(5n+4-\dfrac{5k(k+1)}{2}\right)\right)q^n=: \sum_{n\geq 0} f(n)q^n.\label{eq54}
\end{align}
We apply (\ref{eq31}) and (\ref{eq32}) on (\ref{eq54}), yielding
{\small \begin{align}
\sum_{n\geq 0} f(n)q^n &= f_5f_{10}^2\cdot \dfrac{f_{25}f_{50}^{10}}{f_{10}^{12}}\left(\dfrac{1}{R_5}-q-q^2R_5\right)\nonumber\\
&\cdot\left(\dfrac{1}{R_{10}^4}+\dfrac{q^2}{R_{10}^3}+\dfrac{2q^4}{R_{10}^2}+\dfrac{3q^6}{R_{10}}+5q^8-3q^{10}R_{10}+2q^{12}R_{10}^2-q^{14}R_{10}^3+q^{16}R_{10}^4\right)^2.\label{eq55}
\end{align}}%
We consider the terms of (\ref{eq55}) involving $q^{5n+2}$, divide both sides by $q^2$, and then replace $q^5$ with $q$. In view of (\ref{eq33})--(\ref{eq36}) and the initial values (\ref{eq37})--(\ref{eq310}), we see from (\ref{eq55}) that
\begin{align}
\sum_{n\geq 0} f(5n+2)q^n &= q^3\dfrac{f_1f_5f_{10}^{10}}{f_2^{10}}(-P(3,-2)+2P(3,-1)-10P(2,-1)\nonumber\\
&-16P(1,-1)+27P(1,0)-15)\nonumber\\
&=q^3\dfrac{f_1f_5f_{10}^{10}}{f_2^{10}}\cdot\dfrac{(K-4)^2(K+1)(K^2-3K+1)}{K^2}\nonumber\\
&=q^3\dfrac{f_1f_5f_{10}^{10}}{f_2^{10}}\left(\dfrac{f_1^3f_5}{qf_2f_{10}^3}\right)^2\left(\dfrac{f_2^4f_5^2}{qf_1^2f_{10}^4}\right)\left(\dfrac{qf_1f_{10}^5}{f_2f_5^5}\right)^2(K^2-3K+1)\nonumber\\
&=q^2\dfrac{f_1^7f_{10}^{10}}{f_2^{10}f_5^5}\left((K+1)^2-5(K+1)+5\right)\nonumber\\
&=q^2\dfrac{f_1^7f_{10}^{10}}{f_2^{10}f_5^5}\left(\dfrac{f_2^8f_5^4}{q^2f_1^4f_{10}^8}-5\dfrac{f_2^4f_5^2}{qf_1^2f_{10}^4}+5\right)\nonumber\\
&=\dfrac{f_1^3f_{10}^2}{f_2^2f_5}-5q\dfrac{f_1^5f_{10}^3}{f_2^6f_5^3}+5q^2\dfrac{f_1^7f_{10}^{10}}{f_2^{10}f_5^5}\nonumber\\
&\equiv \dfrac{f_1^3f_{10}^2}{f_2^2f_5}-5q\dfrac{f_{10}^2}{f_2f_5^2}+5q^2\dfrac{f_1^2f_{10}^8}{f_5^4}\pmod{25}.\label{eq56}
\end{align}
Invoking (\ref{eq31}), (\ref{eq32}), and the generating function of $C(n)$, we write (\ref{eq56}) modulo $25$ as
{\footnotesize \begin{align*}
\sum_{n\geq 0} f(5n+2)q^n &\equiv \dfrac{f_{10}^2}{f_5}\sum_{n\geq 0}C(n)q^n+5q^2\dfrac{f_{25}^2f_{10}^8}{f_5^4}\left(\dfrac{1}{R_5}-q-q^2R_5\right)^2\nonumber\\
&-5q\dfrac{f_{50}^5}{f_5^2f_{10}^4}\left(\dfrac{1}{R_{10}^4}+\dfrac{q^2}{R_{10}^3}+\dfrac{2q^4}{R_{10}^2}+\dfrac{3q^6}{R_{10}}+5q^8-3q^{10}R_{10}+2q^{12}R_{10}^2-q^{14}R_{10}^3+q^{16}R_{10}^4\right).
\end{align*}}%
We now look at the terms of the resulting expression involving $q^{5n+4}$. We infer from Theorem \ref{thm12} that
\begin{align*}
\sum_{n\geq 0} f(25n+22)q^n &\equiv \dfrac{f_2^2}{f_1}\cdot 5\dfrac{f_1^2f_5f_{10}^2}{f_2^4}-5\dfrac{f_5^2f_2^8}{f_1^4}-25q\dfrac{f_{10}^5}{f_1^2f_2^4}\\
&\equiv 5f_1^6f_2^8 - 5f_1^6f_2^8\equiv 0\pmod{25}.
\end{align*}
Hence, we arrive at $f(25n+22)\equiv 0\pmod{25}$ for all $n\geq 0$, which proves (\ref{eq12}).
\end{proof}

\section{Proof of Theorem \ref{thm16}}\label{sec5}

\begin{proof}[Proof of Theorem \ref{thm16}] Recall from the proof of (\ref{eq11}) of Theorem \ref{thm15} that $\sum_{n\geq 0} a(5n+1)q^n\equiv 3f_1f_2^2\pmod{5}$. We know from the Jacobi triple product identity \cite[(1.1.1)]{hirsc} that
\begin{align*}
	\dfrac{f_1^2}{f_2} = \sum_{k=-\infty}^\infty (-1)^kq^{k^2}.
\end{align*}	 
Thus, we have 
\begin{align*}
3f_1^3f_2 \equiv 3f_1f_2^2\cdot \dfrac{f_1^2}{f_2}&\equiv\left(\sum_{k=-\infty}^\infty (-1)^kq^{k^2}\right)\left(\sum_{j\geq 0} a(5n+1)q^j\right)\\
&\equiv \sum_{n\geq 0}\left(\sum_{k=-\infty}^\infty (-1)^ka(5n+1-5k^2)\right)q^n\pmod{5}.
\end{align*}
We now write $f_1^3f_2 :=\sum_{n\geq 0} h(n)q^n$ so that
\begin{align}
3h(n) \equiv \sum_{k=-\infty}^\infty (-1)^ka(5n+1-5k^2)\pmod{5}.\label{eq57}
\end{align}
By \cite[Theorem 2]{cophir}, if $p\equiv 13,17,19,23\pmod{24}$, then 
\begin{align*}
h\left(pn+\dfrac{5(p^2-1)}{24}\right) = \pm ph\left(\dfrac{n}{p}\right)
\end{align*}
for all $n\geq 0$, where we denote $h(n/p) = 0$ if $p\nmid n$. Replacing $n$ with $pn+r$ in the above expression, we deduce that
\begin{align}
h\left(p^2n+pr+\dfrac{5(p^2-1)}{24}\right) = 0\label{eq58}
\end{align}
for $r\in \{1,\ldots,p-1\}$. The desired congruence now follows from (\ref{eq57}) and (\ref{eq58}).
\end{proof}

\section{Closing remarks}\label{sec6}
We have shown in this paper new elementary proofs for some of the main results of Amdeberhan and Merca \cite{amdb}, particularly Theorems \ref{thm13} and \ref{thm14}. We have also established congruences modulo $5$ and $25$ for certain finite sums involving $C(n)$ and $a(n)$ via elementary means. Following the proofs of Theorems \ref{thm13} -- \ref{thm16}, we deduce that for $n\geq 0$,
\begin{align*}
\sum_{k =-\infty}^\infty (1+6k)a\left(25n+21-\dfrac{5k(3k+1)}{2}\right)&\equiv 0\pmod{5}.
\end{align*}
We also deduce that for all primes $p\equiv 7,11\pmod{12}, r\in \{1,\ldots, p-1\}$, and $n\geq 0$, 
\begin{align*}
\sum_{k =-\infty}^\infty (-1)^k(3k+1)a\left(5p^2n+5pr+\dfrac{65p^2-41}{24}-5k(3k+2)\right)&\equiv 0\pmod{5}.
\end{align*}
To prove the second identity, we multiply the generating function of $a(5n+1)$ modulo $5$ found in the proof of Theorem \ref{thm15} by the identity \cite[(10.7.7)]{hirsc}
\begin{align*}
	\dfrac{f_2^5}{f_1^2} = \sum_{k=-\infty}^\infty (-1)^k(3k+1)q^{k(3k+2)}
\end{align*}
and then apply the result of Ahlgren \cite{ahl} on the coefficients of the $q$-expansion of $f_2^7/f_1$. We leave the detailed proofs of the above congruences and the discovery of such congruences involving $C(n)$ and $a(n)$ modulo other prime powers to the interested reader.

\section*{Acknowledgment}
The author would like to thank Prof. Jeremy Lovejoy for bringing the paper \cite{choi} to his attention.


\begin{thebibliography}{99}
\bibitem{ahl} S. Ahlgren, {\it Multiplicative relations in powers of Euler's product}, J. Number Theory {\bf 89} (2001), 222--233.	
	
\bibitem{amdb} T. Amdeberhan and M. Merca, {\it From cranks to congruences}, preprint (2025), arXiv:2505.19991.

\bibitem{and} G. E. Andrews, {\it Integer partitions with even parts below odd parts and the mock theta functions}, Ann. Comb {\bf 22} (2018), 433--445.

\bibitem{andg} G. E. Andrews and F. G. Garvan, {\it Dyson’s crank of a partition}, Bull. Amer. Math. Soc. (N.S.) {\bf 18} (1988), 167--171.

\bibitem{atk} A. O. L. Atkin, {\it Proof of a conjecture of Ramanujan}, Glasg. Math. J. {\bf 8} (1967), 14--32.

\bibitem{chehir} S. Chern and M. D. Hirschhorn, {\it Partitions into distinct parts modulo powers of $5$}, Ann. Comb. {\bf 23} (2019), 659--682.

\bibitem{chetn} S. Chern and D. Tang, {\it The Rogers-Ramanujan continued fraction and related eta-quotient representations}, Bull. Aust. Math. Soc. {\bf 103} (2021), 248--259.

\bibitem{choi} D. Choi, S.-Y. Kang and J. Lovejoy, {\it Partitions weighted by the parity of the crank}, J. Combin. Theory Ser. A {\bf 116} (2009), 1034--1046.

\bibitem{cophir} S. Cooper, M. D. Hirschhorn and R. Lewis, {\it Powers of Euler's product and related identities}, Ramanujan J. {\bf 4} (2000), 137--155.

\bibitem{dar} H. B. C. Darling, {\it On Mr. Ramanujan's congruence properties of $p(n)$}, Proc. Cambridge Math. Soc. {\bf 19} (1919), 217--218.

\bibitem{dys} F. J. Dyson, {\it Some guesses in the theory of partitions}, Eureka (Cambridge) {\bf 8} (1944), 10--15.

\bibitem{garv} F. G. Garvan, {\it Generalizations of Dyson's rank}, Ph. D. thesis, Pennsylvania State University, 1986. 

\bibitem{hirsc} M. D. Hirschhorn, {\it The Power of $q$, A Personal Journey}, Developments of Mathematics, vol. 49, Springer, Cham, 2017.

\bibitem{radu} C.-S. Radu, {\it An algorithmic approach to Ramanujan–Kolberg identities}, J. Symbolic Comput. {\bf 68} (2015), 225--253.

\bibitem{rama} S. Ramanujan, {\it Some properties of $p(n)$, the number of partitions of $n$}, Proc. Cambridge Math. Soc. {\bf 19} (1919), 210--213.

\bibitem{smoot} N. A. Smoot, {\it On the computation of identities relating partition numbers in arithmetic progressions with eta quotients: An implementation of Radu's algorithm}, J. Symbolic Comput. {\bf 104} (2021), 276--311.

\bibitem{tang} D. Tang, {\it Congruences modulo powers of 5 for the crank parity function}, Quaestiones Math. {\bf 47} (2024), 715--734.

\bibitem{wat} G. N. Watson, {\it Ramanujans Vermutung \"{u}ber Zerf\"{a}llungsanzahlen}, J. Reine Angew. Math. {\bf 179} (1938), 97--128.
\end{thebibliography}
\end{document}